\documentclass[11pt,a4paper]{article}

\usepackage{latexsym}
\usepackage{amssymb}
\usepackage[mathscr]{eucal}
\usepackage{epsfig}
\usepackage{lscape}
\usepackage{amssymb}
\input amssym.def
\input amssym

\newtheorem{Theorem}{Theorem}[section]

\newtheorem{Remark}[Theorem]{Remark}
\def\proof#1. {\par
                      \ifdim\lastskip<15pt
                      \removelastskip\penalty-200
                      \vskip5pt plus3pt minus3pt
                      \fi
                       {\def\a{#1}
                       \ifx\a\empty
                       {\noindent\bf Proof.}
                       \else
                       {\noindent\bf Proof of #1.}
                       \fi}\enspace}

\def\endproof{\hfill\hspace{-6pt}\rule[-4pt]{6pt}{6pt}
\vskip8pt plus3pt minus 3pt}

\title{Convexity of the zeros of some orthogonal polynomials and related functions}

\author{Kerstin Jordaan\thanks{Department of
Mathematics and Applied Mathematics, University of Pretoria, Pretoria, 0002, South Africa.
Research by this author is partially supported by the National Research Foundation under grant
number 2054423.} \and Ferenc Tookos\thanks{Institute for Biomathematics and Biometry, GSF -
Research Center for Environment and Health, Neuherberg, Germany.
Research by this author is partially supported by OTKA 49448.} }
\date{}
\begin{document}

\maketitle
\smallskip

\begin{center}
\end{center}
\medskip

\begin{abstract} We study convexity properties of the zeros of some
special functions that follow from the convexity theorem of Sturm.
We prove results on the intervals of convexity for the zeros of
Laguerre, Jacobi and ultraspherical polynomials, as well as
functions related to them, using transformations under which the
zeros remain unchanged. We give upper as well as lower bounds for
the distance between consecutive zeros in several cases.
\end{abstract}
\bigskip

\noindent AMS MOS Classification:\quad 33C45, 34C10, 42C05

\newpage
\section{Introduction}
The Sturm comparison theorem for solutions of second order
differential equations of the form $y''+F(t)y=0$ (cf. \cite{Sturm})
has been significantly extended since its publication 170 years ago.
Some of the immediate applications to zeros of the solutions $y(t)$,
and those of the derivative $y'(t)$, include Sonin's theorem on the
monotonicity of extrema of such solutions (cf. \cite{Makai}), and a
result known as Sturm's convexity theorem, first mentioned in
\cite{Sturm}, on the monotonicity of distances between the zeros of
the solution (cf. \cite{Laforgia}, \cite{Makai} and \cite{Szego2}).

\medskip Sonin's theorem was extended to more general differential equations of the form
$P(t)y''+Q(t)y'+y=0$ using a remarkably simple proof (cf.
\cite{Redheffer}, p.443). In this form the theorem can be directly
applied to classical orthogonal polynomials such as Hermite,
Laguerre and Jacobi polynomials, providing the monotonicity of their
relative maximum values and estimates on their supremum norm. This
has been done from a different perspective in \cite{Kalamajska},
recovering results for Legendre, Laguerre and Jacobi polynomials
given in \cite{Szego}.

\medskip
In this paper, we consider the implications of the convexity theorem
of Sturm for the convexity of the zeros and bounds on the distance
between the zeros of some classical orthogonal polynomials and
functions related to them.





\section{Convexity and spacing of zeros}
The convexity theorem and an obvious consequence of the comparison
theorem of Sturm, already noted in \cite{Sturm}, can be summarised
as follows.

\begin{Theorem}\cite{D-G-S}\label{convexity}  Let $y''(t)+F(t)y(t)=0$ be a second-order differential equation in normal form,
where $F$ is continuous in $(a,\ b)$. Let $y(t)$ be a nontrivial
solution in $(a,\ b)$, and let $x_1<\ldots<x_k<x_{k+1}<\ldots$
denote the consecutive zeros of $y(t)$ in $(a,\ b)$. Then
\begin{enumerate}
\item if $F(t)$ is strictly increasing in $(a,\ b)$, $x_{k+2}-x_{k+1}<x_{k+1}-x_k$,
\item if $F(t)$ is strictly decreasing in $(a,\ b)$, $x_{k+2}-x_{k+1}>x_{k+1}-x_k$.
\item if there exists $M>0$ such that $F(t)<M$ in $(a,\ b )$ then
\[\Delta x_k\equiv x_{k+1}-x_k>\frac{\pi}{\sqrt{M}},\]
\item if there exists $m>0$ such that $F(t)>m$ in $(a,\ b )$ then
\[\Delta x_k<\frac{\pi}{\sqrt{m}},\]
\end{enumerate}
\end{Theorem}
We say that the zeros of $y$ are concave (convex) on $(a,b)$ for the
first (second) case.


\medskip The convexity theorem has been used to obtain the variation of convexity
properties with respect to a parameter, or the order, for the zeros
of gamma, q-gamma, Bessel, cylindrical and Hermite functions as
described in the survey paper \cite{Muldoon}.


\medskip
In order to apply the convexity theorem to special functions that
are solutions of second-order differential equations, the
differential equation has to be transformed into normal form. One
simple way to do this is through the following change of dependent
variable. Let

\[x''+g(t)x'+f(t)x=0\]
be a second-order differential equation and set

\begin{equation} y=x\exp\left( \frac12\int^t g(s)ds\right) \label{transform}\end{equation}

The corresponding equation for $y$ is in normal form:

\[y''+F(t)y=0,\]
where $F(t)=f(t)-\frac14 g^2(t)-\frac12 g'(t)$. The advantage of
this transformation is that it does not change the independent
variable, and the zeros of $x$ and $y$ are the same. Hille
\cite{Hille} already used transformation (\ref{transform}) to prove
the convexity of zeros of the Hermite polynomials.

\medskip
It is also possible to consider other changes of variable and obtain
information on the convexity of the transformed zeros. This was done
already by Szeg\H{o} for the ultraspherical polynomials
\cite[Theorem 6.3.3]{Szego} and lately by Dea\~no, Gil and Segura
\cite{D-G-S, D-S} for hypergeometric functions.

\medskip
We will consider the convexity and spacing of the zeros of special
functions such as Laguerre, Jacobi and, as a special case, the
ultraspherical polynomials, for fixed order $n$, by transforming
their differential equations to normal form using (\ref{transform}).
Sturm \cite{Sturm} used the same method to obtain results on the
convexity and spacing of the zeros of the Bessel function.
Interesting work on the spacing of the zeros of Jacobi polynomials,
as the degree changes, is done in \cite{Laforgia2}. For higher
monotonicity refer to, amongst others, \cite{Lorch} and
\cite{Bessel}.

\medskip We note that since the convexity theorem is applicable to any
oscillating solutions of second order differential equations in
normal form, the results we obtain are not restricted to the
polynomial cases, i.e. $n$ need not necessarily be an integer, as
long as the corresponding functions are oscillating on the interval
under consideration. In addition, the results can be extended to
parameter values where the polynomials are no longer orthogonal,
since quasi-orthogonality ensures the existence of some zeros on the
interval of orthogonality (cf. \cite{Brez}).


\section{Laguerre polynomials}
The differential equation
\[tx''+(\alpha+1-t)x'+nx=0\] satisfied by the Laguerre polynomials, $L_n^{\alpha}(t)$, orthogonal on
$(0,\infty)$ with respect to the weight function $t^{\alpha}e^{-t}$
when $\alpha>-1$, is transformed to
\begin{eqnarray}y''+F(t)y&=&0\nonumber\end{eqnarray} by (\ref{transform})
where \begin{eqnarray}F(t)&=&\frac{-t^2+2\alpha
t+2t+4nt-\alpha^2+1}{4t^2}\label{F_Laguerre}.\end{eqnarray} $F'(t)$
changes sign at
\begin{equation}
t_0:=\frac{\alpha^2-1}{\alpha+2n+1}.\label{t_0}
\end{equation}
\begin{Theorem}\label{L}The zeros of $L_n^{\alpha}(t)$ on $(0,\
\infty)$ are \begin{enumerate} \item all convex if $n>0$ and
$-1<\alpha\leq 3$ \item all convex if $\alpha>3$ and
$0<n<\frac{\alpha+1}{\alpha-3}$ \item concave for $t<t_0$ and convex
for $t>t_0$ when $\alpha>3$, $n>\frac{\alpha+1}{\alpha-3}$ and $t_0$
is defined by (\ref{t_0}) .
\end{enumerate}
Moreover, for the distance between consecutive zeros we have the
general estimate
\begin{equation} \Delta x_k>\frac{\pi\sqrt{2}}{\sqrt{2\alpha n+\alpha+2n^2+2n+1}}
\qquad k=1,\ldots, n-1\label{Lag-gen}\end{equation}
and also if $x_k>t_0$ then
\begin{equation}\Delta x_k>\frac{\pi}{\sqrt{F(x_k)}} \qquad k=1,\ldots, n-1\label{Lag-low}\end{equation}
and
\begin{equation}\Delta x_k<\frac{\pi}{\sqrt{F(x_{k+1})}} \qquad k=1,\ldots, n-2\label{Lag-upp}\end{equation}
where $F$ is defined by (\ref{F_Laguerre}).
\end{Theorem}
\proof. For $|\alpha|<1$, $t_0<0$, hence $F(t)$ will be decreasing
on ($0,\infty)$. When $\alpha\geq1$, $F(t)$ is increasing on $(0,\
t_0)$ and decreasing on $(t_0,\ \infty)$.
Let $x_1$ denote the smallest zero of $L_n^{\alpha}$, then we know
that $x_{1}>\frac{\alpha+1}n$ (cf. \cite{Hahn}). This implies that
when $t_0<\frac{\alpha+1}n$, $F(t)$ will be decreasing on the
interval $(x_1,\infty)$. An easy calculation shows that this
condition is equivalent to either $\alpha\leq 3$ or $\alpha>3$ and
$n<\frac{\alpha+1}{\alpha-3}$. The estimates on the distance $\Delta
x_k$ follow from Theorem \ref{convexity}(3),(4). The maximum of $F$
is at $t_0$ and $F(t_0)>0$, therefore we can take $F(t_0)$ as $M$ to
obtain (\ref{Lag-gen}). For (\ref{Lag-low}) and (\ref{Lag-upp}), we
use the fact that when $x_k>t_0$, $F$ is monotone decreasing on
$(x_k,\ x_{k+1})$. In fact, $F$ is monotone decreasing on
$(0,\infty)$ and tends to -1/4 as $t\to\infty$, so there is exactly
one point $t_1$ on $(t_0,\ \infty)$, where $F$ crosses the $x$-axis.
The form of the differential equation implies that if $F(t)<0$ and
$y(t)>0$, the graph will be concave up and similarly, if $y(t)<0$,
the graph will be concave down. Hence there can be at most one zero
of the Laguerre polynomial to the right of $t_1$. This means that
$F(x_{n-1})$ is positive, but $F(x_n)$ may be negative and therefore
the index in (\ref{Lag-upp}) only runs up to $n-2$.\endproof

\begin{Remark} An interesting question is whether it is possible to
find $\alpha$ and $n$ values so that the first several zeros of the
Laguerre polynomial are concave. This would require $t_0$ to be
greater than $x_{1}$. In this regard we note that $t_0$ is always
less than $\frac{(\alpha+1)(\alpha+2)}{\alpha+n+1}$, which is the
upper bound for $x_{1}$ given in \cite{Hahn}. It is also always less
than $\frac{(\alpha+1)(\alpha+3)}{\alpha+2n+1}$, the upper bound
given in \cite{Szego}.
\end{Remark}

\section{Jacobi polynomials}
The differential equation for Jacobi polynomials, $P_n^{(\alpha,
\beta)}(t)$, orthogonal on $(-1,1)$ with respect to the weight
function $(1-t)^{\alpha}(1+t)^{\beta}$ when $\alpha,~\beta>-1$, is
\[(1-t^2)x''(t)+(\beta-\alpha-(\alpha+\beta+2)t)x'(t)+n(n+\alpha+\beta+1)x(t)=0.\]
In the normal form, $y''+F(t)y=0$, we have
\begin{eqnarray}
F(t)&=&\frac{-zt^2-2(x-y)t-2x-2y+z}{4(t^2-1)^2}\label{F(t)}\\
\mbox{with}~x&=&\alpha^2-1\nonumber\\
y&=&\beta^2-1\nonumber\\
z&=&(\alpha+\beta+2n)(\alpha+\beta+2n+2).\nonumber\end{eqnarray}
Also \begin{eqnarray*}
F'(t)=\frac{zt^3+3(x-y)t^2+(4x+4y-z)t+(x-y)}{2(t^2-1)^3}:=\frac{j(t)}{2(t^2-1)^3}
\end{eqnarray*}
and we denote the discriminant of $j'(t)$ by
\[D:=12(3x^2+3y^2+z^2-6xy-4xz-4yz).\]
For the convexity theorem to be applicable, we need oscillating
solutions. The condition on the parameters for this is (cf.
\cite{D-G-S}) \[n>0,\ n+\alpha+\beta>0,\ n+\alpha>0,\ n+\beta>0.\]
From now on we shall assume that the coefficients satisfy these
conditions.
\begin{Theorem}\label{Jacobi}
If $|\alpha|>1,~|\beta|<1$ and $D<0$, all the zeros of
$P_n^{(\alpha,\beta)}$ on the interval $(-1,1)$ are
convex.\end{Theorem}

\proof. $F(t)$ is a rational function with vertical asymptotes at
$t=\pm 1$. If $|\alpha|>1$ and $|\beta|<1$, then $j(-1)=-8y>0$ and
$j(1)=8x>0$, so that $\lim_{t\to -1}F(t)=\infty$ and $\lim_{t \to
1}F(t)=-\infty$. $D<0$ implies that $j'(t)\neq 0$ for $t \in(-1,1)$
and hence $j(t)$ will have no extreme values on this interval. It
follows that $F(t)$ is monotone decreasing on $(-1,1)$ and Theorem
\ref{convexity}(2) yields the result.
\endproof
Note that the conditions of Theorem \ref{Jacobi} are satisfied if,
for example, $y$ is sufficiently small and $x<z<3x$. This is true
if, for instance, $\beta$ is sufficiently close to $-1$, $\alpha>1$,
and $n<\frac12 (-\alpha+\sqrt{3\alpha^2-2})$. Also $D<0$ will be
satisfied for sufficiently large $\alpha$ if we fix $\beta$ and $n$.
However, for fixed $\alpha$ and $\beta$ the discriminant $D$ is
positive for large $n$.

\medskip More results of this type can be obtained by
ensuring the positivity of $j(t)$ on $(-1,1)$. Denote the zeros of
$j'(t)$ by \[t_{1,2}=\frac{6(y-x)\pm \sqrt{D}}{6z},\] then we have
that the zeros of $P_n^{(\alpha,\beta)}$ on $(-1,\ 1)$ are convex if
$|\alpha|>1$ and $|\beta|<1$ and if $t_i\in (-1,1)$ for some $i=1,2$
then $j(t_i)>0$. One can prove conditions for concavity of the zeros
of $P_n^{(\alpha,\beta)}$ on the whole interval $(-1,\ 1)$ in a
similar manner.

\medskip The general study of the convexity of the zeros is difficult,
since the convexity intervals are determined by the roots of $j(t)$,
which are hard to handle due to the 3 parameters. However, there are
still some things that can be said about the general case.

\medskip
The oscillation condition $z>0$ implies that $j$ has a concave part,
followed by a convex, the inflection point being at
$t_0:=\frac{y-x}{z}$. Therefore it follows from the shape of $j(t)$
that there may be up to 4 different intervals of changing concavity
for the zeros on $(-1,\ 1)$: (from left to right)
concave-convex-concave-convex. Any of these intervals may be missing
from the sequence, for example, concave-convex-concave or
convex-concave are possible for certain parameter values.

\medskip It is interesting to analyse the convexity of the zeros of
$P_n^{(\alpha,\beta)}$ for sufficiently
large degree.
\begin{Theorem} Let $\alpha$ and $\beta$ be fixed and let $n \to
\infty$, then the convexity of the zeros of $P_n^{(\alpha,\beta)}$
on $(-1,\ 1)$ changes in the following way (from left to right):
\begin{enumerate}
\item if $|\alpha|\leq 1$ and $|\beta|\leq 1$ then convex-concave,
\item if $|\alpha|\leq 1$ and $|\beta|>1$ then concave-convex-concave,
\item if $|\alpha|>1$ and $|\beta|\leq 1$ then convex-concave-convex,
\item if $|\alpha|>1$ and $|\beta|>1$ then concave-convex-concave-convex.
\end{enumerate}
\end{Theorem}
\proof. If $\alpha$ and $\beta$ are fixed and $n\rightarrow\infty$,
an easy calculation shows that the local extremum locations of
$j(t)$ tend to $\pm \sqrt{3}/3$. Since $z>0$, the local extremum
near $t=-\sqrt{3}/3$ will be the maximum and this maximum value
tends to $\infty$. Similarly, the minimum value near $\sqrt{3}/3$
tends to $-\infty$. Since the inflection point
$t_0=\frac{y-x}{z}\rightarrow 0$, there is at least one change of
concavity in $(-1,\ 1)$ (from convex to concave) and whether there
are more, depends on the sign of $j(-1)$ and $j(1)$. Now
$j(-1)=8(1-\beta^2)$ and $j(1)=8(\alpha^2-1)$ and the result
follows.
\endproof
\section{Ultraspherical polynomials} An important special case of
the Jacobi polynomials are the ultraspherical polynomials,
$P_n^{(\alpha,\alpha)}(t)$ where $\alpha=\beta.$ In this case
\[F(t)=\frac{-(\alpha+n)(\alpha+n+1)t^2+(1+n+n^2+\alpha+2\alpha
n)}{(t^2-1)^2},\]
the numerator of $F'(t)$ is
$j(t)=4[(\alpha+n)(\alpha+n+1)t^3-(2+n+n^2+\alpha+2\alpha
n-\alpha^2)t]$
and the discriminant of $j'(t)$ is
$D=192(\alpha+n)(\alpha+n+1)(2+n+n^2+\alpha+2\alpha n-\alpha^2).$
Note that the leading coefficient of $j(t)$ is positive when
$\alpha>-1$. The point of inflection of $j(t)$ is $t_0=0$ and hence
the convexity of zeros changes exactly in the middle of the interval
$(-1,1)$. The local extrema of $j(t)$ are at
\begin{eqnarray*}t_{1,2}&=&
\pm\sqrt{\frac{(n+\alpha)(n+\alpha+1)-2(\alpha^2-1)}{3(n+\alpha)(n+\alpha+1)}},
\end{eqnarray*}
and the two remaining zeros of $j(t)$ are
\begin{equation} T_{1,2}=\pm\sqrt{\frac{(n+\alpha)(n+\alpha+1)-2(\alpha^2-1)}{(n+\alpha)(n+\alpha+1)}}=
t_{1,2}\sqrt{3}.\label{T}\end{equation} where $T_1$ denotes the
negative and $T_2$ the positive zero.
\begin{Theorem} If $|\alpha|\leq 1$, the zeros of $P_n^{(\alpha,\alpha)}$ on $(-1,\ 0)$ are
convex and those on $(0,\ 1)$ are concave. In addition \[\Delta
x_k<\frac{\pi}{\sqrt{F(0)}}=\frac{\pi}{\sqrt{2\alpha
n+\alpha+n^2+n+1}},\] and for the positive zeros we have
\[\frac{\pi}{\sqrt{F(x_{k+1})}}<\Delta x_k<\frac{\pi}{\sqrt{F(x_{k})}}.\]
\end{Theorem}
\proof. For $|\alpha|\leq 1$ we have $T_1<-1$ and $T_2>1$, so $j$ is
positive on $(-1,\ 0)$ and negative on $(0,\ 1)$. Therefore $F(t)$
is decreasing on $(-1,0)$ and increasing on $(0,1)$ and the
convexity of the zeros follows from Theorem \ref{convexity}(1),(2).
In addition, $F(0)>0$ is a minimum value, so we have an upper bound
on the distance between any two consecutive zeros from Theorem
2.1(4). Finally, since $F(t)$ is increasing on $(0,1)$,
$0<F(x_i)<F(x)<F(x_{i+1})$ for $x\in(x_i,x_{i+1})$, where
$x_i~x_{i+1}$ are any two consecutive positive zeros and the last
inequality follows from Theorem \ref{convexity}(3),(4).
\endproof
Note that this result also applies to Chebyshev and Legendre
polynomials as special cases with $\alpha=\frac 12$ and $\alpha=0$
respectively.
\begin{Theorem}
Let $|\alpha|>1$ and $(n+\alpha)(n+\alpha+1)\leq 2(\alpha^2-1)$ then
the zeros of $P_n^{(\alpha,\alpha)}$ on $(-1,\ 0)$ are concave and
those on $(0,\ 1)$ are convex. Furthermore
\[\Delta x_k>\frac{\pi}{\sqrt{F(0)}}.\]
\end{Theorem}
\proof. If $|\alpha|>1$, and $D<0$, $j$ has no local extremum and is
monotone increasing on $(-1,1)$. This, together with the fact that
$j(0)=0$, implies that $F(t)$ is increasing on $(-1,0)$ and
decreasing on $(0,1)$, and the result follows. \endproof

Note that we cannot obtain estimates involving $x_k,\ x_{k+1}$,
because $F\rightarrow -\infty$ as $t \to \pm 1$.


\begin{Theorem}
Let $|\alpha|>1$ and $(n+\alpha)(n+\alpha+1)>2(\alpha^2-1)$ then the
zeros of $P_n^{(\alpha,\alpha)}$ are concave on $(-1,\ T_1)$ and
$(0,\ T_2)$ and convex on $(T_1,\ 0)$ and $(T_2,\ 1)$, where
$T_{1,2}$ are as in (\ref{T}). We also have that
\[\Delta x_k>\frac{\pi}{\sqrt{F(T_2)}},\]
moreover, if $(x_k,\ x_{k+1})\subset(T_1,\ T_2)$ then
\[\Delta x_k<\frac{\pi}{\sqrt{F(0)}},\]
and if $(x_k,\ x_{k+1})\subset(0,\ T_2),$ then
\[\frac{\pi}{\sqrt{F(x_{k+1})}}<\Delta x_k<\frac{\pi}{\sqrt{F(x_{k})}}.\]
\end{Theorem}
\proof. When $|\alpha|>1$ and $D>0$, $j(t)$ has 3 zeros on $(-1,1)$,
namely at $T_1,~0$ and $T_2$ with $F(t)$ having local maxima at
$T_{1,2}$ and a local minimum at $t=0$.
\endproof

\end{document}